\documentclass[12pt,a4paper]{article}
\usepackage{amstext}
\usepackage{fancyhdr}
\usepackage{amsfonts,graphicx,bezier, amssymb}
\usepackage{amsmath}
\usepackage{caption}
\usepackage{mathtools}
\usepackage{fancyhdr}
\usepackage{float}
\usepackage{pgfplots}
\pgfplotsset{compat=newest}
\usepackage[utf8]{inputenc}
\newcommand{\noi}{\noindent}
\newtheorem{thm}{Theorem}[section]
\newtheorem{defn}{Definition}[section]

\newtheorem{rem}{Remark}[section]
\newtheorem{lm}{Lemma}[section]
\newtheorem{ex}{Example}[section]
\newtheorem{cor}{Corollary}[section]
\newtheorem{prop}{Proposition}[section]
\newenvironment{prf}{\noindent{\bf{Proof:}}~~}{\hfill\rule{1ex}{1ex}\vskip1.5ex}

\usepackage{tikz}
\usepackage{color}
\usepackage{xcolor}
\usetikzlibrary{arrows}
\usetikzlibrary{matrix}

\usepackage{etoolbox}
\usepackage[colorlinks]{hyperref}

\begin{document}
	
	\tikzset{ 
		table/.style={
			matrix of nodes,
			row sep=-\pgflinewidth,
			column sep=-\pgflinewidth,
			nodes={
				rectangle,
				draw=black,
				align=center
			},
			minimum height=0.5em,
			text depth=4.0ex,
			text height=2ex,
			nodes in empty cells,
			
			every even row/.style={
				nodes={fill=gray!20}
			},
			column 1/.style={
				nodes={text width=1em,font=\bfseries}
			},
			row 1/.style={
				nodes={
					fill=white,
					text=black,
					font=\bfseries
				}
			}
		}
	}

		
		\begin{center}\Large{\bf{Reduced submodules of finite dimensional polynomial modules}}
			
		\end{center}
		\vspace*{0.3cm}
		\begin{center}
			Tilahun Abebaw, Nega Arega, Teklemichael Worku Bihonegn \footnote{This work forms part of the third author's PhD thesis.} and David Ssevviiri \footnote{Corresponding author.}
		\end{center}
		\begin{abstract}
			Let $k$ be a field with characteristic zero, $R$ be the ring $k[x_1, \cdots, x_n]$ and $I$ be a monomial ideal of $R$. We study the Artinian local algebra $R/I$ when considered as an $R$-module $M$. We show that the largest reduced submodule of $M$ coincides with both the socle of $M$ and the $k$-submodule of $M$ generated by all outside corner elements of the Young diagram associated with $M$. Interpretations of different reduced modules is given in terms of Macaulay inverse systems. It is further shown that these reduced submodules are examples of modules in a torsion-torsionfree class, together with their duals; coreduced modules, exhibit symmetries in regard to Matlis duality and torsion theories. Lastly, we show that any $R$-module $M$ of the kind described here satisfies the radical formula. 
			
		\end{abstract}
		
		{\bf Keywords}: Reduced modules, socle, Macaulay inverse systems, local Artinian algebras, modules that satisfy the radical formula  
		
		\vspace*{0.4cm}
		
		\section{Introduction}
		\begin{paragraph}\noi
			
			Reduced rings play an important role in algebra. A module analogue of reduced rings was defined by Lee and Zhou in \cite{lee2004reduced}. Reduced modules have since been studied in \cite{kyomuhangi2021generalised,kyomuhangi2020locally,rege2008reduced, ssevviiri2022applications,ssevviiri2023applications} among others.
			Let $R$ be a commutative unital ring and $I$ be an ideal of $R$. Reduced modules form a full subcategory of $R$-Mod on which the $I$-torsion functor $\Gamma_I$ is representable, i.e., if $M$ is an $I$-reduced $R$-module, then $\Gamma_I(M)\cong \text{Hom}(R/I, M)$, see \cite{ssevviiri2022applications}. It was shown in \cite{ssevviiri2022applications} that reduced modules and their dual, coreduced modules, provide a setting in which both the Matlis-Greenlees-May Equivalence and Greenlees-May Duality hold. In \cite{kimuli2022characterizations}, reduced modules were used to characterize regular modules. In \cite{ssevviiri2023applications}, it was shown that $I$-reduced and $I$-coreduced modules provide the necessary and sufficient conditions for the functor $\text {Hom}_R(R/I, -)$ to be a radical. The same conditions unify and subsume different conditions which were proved on a case-by-case basis for the $I$-torsion functor $\Gamma_I$ to be a radical. A more general version of reduced modules was studied in \cite{kyomuhangi2021generalised}. We note that this is what Rohrer and Yekutieli studied in \cite{rohrer2019torsion} and \cite{yekutieli2021weak} respectively, although called them modules with bounded torsion. The same modules are called modules whose submodules $(0:_M a^t)$ with $a \in R$ and $t\in \mathbb{Z}^+$ are stationary, by Schenzel and Simon in \cite[Proposition $3.1.10$]{schenzel2018completion}. This general version of reduced modules relates to prisms which belong to the groundbreaking theory of perfectoid rings, \cite{bhatt2019prisms,yekutieli2021weak}.
		\end{paragraph}
		\begin{paragraph}\noi
			A module $M$ over a commutative ring $R$ is {\it reduced} if for all $a\in R$ and $m\in M,~ a^2m=0$ implies that $am=0$.
			We wish to study for an arbitrary $R$-module $M$ in a suitable full subcategory of the category of $R$-modules, the set
			$$
			\mathfrak{R}(M):=\{m\in M: a^2m=0 \Rightarrow am=0 ~\text{for~~all}~ a\in R\}.
			$$ By definition, $M$ is a reduced $R$-module if and only if $\mathfrak{R}(M)=M$. It is easy to see that in general, $\mathfrak{R}(M)$ is not a submodule of $M$, see Example \ref{Rikard}. 
			Let $k$ be a field, $R:=k[x_1,\cdots, x_n]$ and $I$ a monomial ideal of $R$. In this paper, we characterize reduced submodules, $\mathfrak{R}(M)$ of $R$-modules $M$ in the full subcategory $\mathfrak{C}$ of $R$-Mod consisting of $R$-modules of the form $R/I$, with $\text{dim}_k  (R/I)  < \infty $.
			A notion which has been so useful in characterizing $\mathfrak{R}(M)$ is that of socle of $M$. For an Artinian local algebra, $M:=R/I$ and a maximal ideal ${\bf m}$ of $R$, the socle of $M$, $\text{Soc}(M)$ is the submodule of $M$ given by $(0:_M{\bf m})$; the collection of all elements of $M$ annihilated by {\bf m}. This definition is equivalent to saying that $\text{Soc}(M)$ is the direct sum of all simple submodules of $M$. Socle of Artinian local algebras has been widely studied, see \cite{agnarsson2020monomial,agnarsson2023posets,bruns1998cohen,eizenbud1995commutative,van2010simplicial,villarreal7300monomial} among others. It is well known, for instance that a local Artinian algebra $R/I$ is Gorenstein if and only if $\text{dim}_k(\text{Soc}(R/I))=1$.
		\end{paragraph}
		\begin{paragraph}\noi
			We show in Section $2$ that
			for any $M\in \mathfrak{C}$, $\mathfrak{R}(M)$
			is a submodule of $M$ which coincides with both $\text{Soc}(M)$, and with the submodule of $M$ generated by all outside corner elements of the Young diagram associated with $M$, Theorem \ref{thm soc(M)=R(M)}. We hope that this coincidence will increase the versatility of both the largest reduced submodule of $M$ and the socle of $M$, which are already widely studied. In Section $3$, we exploit the Macaulay inverse systems to give a correspondence between different reduced submodules of $M$ in $\mathfrak{C}$ and their associated  Macaulay inverse duals. The correspondences are summarized in Figure \ref{Fig table}. In Section $4$, we exhibit using a diagram, see Figure \ref{fig summary}, symmetries from the following notions: $I$-reduced, $I$-coreduced, $I$-torsion and $I$-complete together with their Matlis duals and associated torsion theories. The work in this section is mainly complementary to that done in papers \cite{ssevviiri2022applications,ssevviiri2023applications}. Lastly, in Section $5$, we show that modules $M$ in the subcategory $\mathfrak{C}$ satisfy the radical formula and their semiprime radicals coincide with the Jacobson radical.
		\end{paragraph}	
		
		\section{For any $M\in \mathfrak{C}, \mathfrak{R}(M) ~\text{coincides~ with Soc}(M)$}
		\begin{paragraph}\noi
			In this section, we prove the coincidence of $\mathfrak{R}(M)$ with both $\text{Soc}(M)$ and the submodule of $M$ generated by all outside corner elements of $M$, for any $M\in \mathfrak{C}$.
		\end{paragraph}
		\begin{defn}\normalfont\label{defn ouer corner }
			Let $R:=k[x_1, \cdots, x_n]$ and $M\in \mathfrak{C}$. An element $m\in M $ is an {\it outside corner element} if $x_im=0$ for all $1\leq i\leq n$. $m\in M$ is {\it inner} if it is not an outside corner element.
		\end{defn}
		\begin{paragraph}\noi
			In Definition \ref{defn ouer corner }, if $n=2$, then this definition is exactly what appears in \cite[page 8]{fulton1997young} and given combinatorially as the boxes which are at the outside corners of the Young diagram. Note that any $M\in \mathfrak{C}$ can be viewed as a $k$-vector space or an $R$-module. When we say, ``generating set of $M$", we will always mean the monomial basis of $M$ seen as a $k$-vector space.
		\end{paragraph}
		
		\begin{lm}\normalfont\label{lem simple implies reduced}
			Every simple module is reduced.
		\end{lm}
		\begin{prf}
			We prove first that a simple module is prime. An $R$-module $M$ is prime if for any $a\in R$ and $m\in M, am=0$ implies that either $m=0$ or $aM=0$. Now, suppose that $M$ is simple, $a\in R$ and $m\in M$ such that $am=0$. Then $aRm=0$. $M$ being simple implies that either $Rm=0$ or $Rm=M$. If $Rm=0$, then $m=0$. Suppose $Rm=M$, then $aM=0$. So, $M$ is prime. We now prove that a prime module is reduced. Let $M$ be a prime $R$-module, $a\in R$ and $m\in M$ such that $a^2m=0$. Then we have $a(am)=0$. By definition of a prime module, we have $am=0$ or $aM=0$. In both cases $am=0$, since $am\in aM$. This shows that $M$ is reduced.
		\end{prf}
		\begin{thm}\label{thm soc(M)=R(M)}
			\normalfont
			Let $M\in \mathfrak{C}$, $S:=\{m_1, m_2,\cdots, m_n\}$ be the collection of all outside corner elements of $M$ and $\langle S\rangle_k$ be the $k$-submodule of $M$ generated by $S$. $$\langle S\rangle_k=\mathfrak{R}(M)=\text{Soc}(M).$$
		\end{thm}
		\begin{prf}
			Let $m\in \langle S \rangle_k$, i.e., $m=\sum_{i=1}^{n}r_im_i$ where $r_i \in k$ and $m_i \in S$. Suppose that $a^2m=0$ for some $a\in R$.
			If $a\in \langle x_1, \cdots, x_n \rangle$, then by definition of outside corner elements, $am=0$. Now suppose that $a\in R\setminus \langle x_1, \cdots, x_n \rangle$. If $a=a_0+a_1$, where $a_0\in k$ and $a_1\in\langle x_1, \cdots, x_n \rangle$, then $am=a_0m\neq 0$. So, $a^2m=a_0^2m\neq 0$. Thus, in all cases $m\in \mathfrak{R}(M)$ and $\langle S \rangle_k \subseteq \mathfrak{R}(M)$. 
			Suppose that $\langle S\rangle_k \subsetneq \mathfrak{R}(M)$. Then  there exists $m\in \mathfrak{R}(M)$ and $m\notin \langle S\rangle_k$. This implies that $m$ is not an outside corner element of $M$. So, there exists $x_i$ for some $i\in \{1,\cdots, n\}$ such that $x_im\neq 0$. However, since $a\in \langle x_1, \cdots, x_n\rangle$ for sufficiently large $t\in \mathbb{Z}^+, a^tm=0$. This shows that $m\notin \mathfrak{R}(M)$ which is a contradiction. It is therefore impossible for the inclusion, $\langle S \rangle_k\subseteq \mathfrak{R}(M)$ to be strict. Thus $\langle S\rangle_k=\mathfrak{R}(M)$. Any simple module is reduced, Lemma \ref{lem simple implies reduced} and a direct sum of reduced modules is reduced, \cite{lee2004reduced}. So, $\text{Soc}(M)\subseteq \mathfrak{R}(M)$. Let $m\in \mathfrak{R}(M)$. Since $\mathfrak{R}(M)=\langle S\rangle_k, m\in \langle S\rangle_k$ and therefore $x_im=0$ for each $1\leq i\leq n$. So, $\langle x_1, \cdots, x_n\rangle m=0$, by definition of socle it follows that $m\in \text{Soc}(M)$. Thus, $\mathfrak{R}(M)\subseteq \text{Soc}(M)$.
		\end{prf}
		\begin{cor}\normalfont
			If $\mathfrak{C}_\text{red}$ is the collection of all reduced $R$-modules $N$ such that $N$ is a submodule of $M\in \mathfrak{C}$, then $\mathfrak{C}_\text{red}$ coincides with all semisimple submodules $N$ of $M\in \mathfrak{C}$.
		\end{cor}
		\begin{prf}
			By Theorem \ref{thm soc(M)=R(M)}, for any $M\in \mathfrak{C}$, $\mathfrak{R}(M)=\text{Soc}(M)$. If $N$ is a submodule of $M\in \mathfrak{C}$ which is reduced, then $N\subseteq \mathfrak{R}(M)$ the largest reduced submodule of $M$. So, $N$ is semisimple. Similarly, if $N$ is a semisimple submodule of $M\in \mathfrak{C}$, then $N\subseteq \text{Soc}(M)$ and therefore, $N$ is a reduced submodule of $M$.
		\end{prf}
		\begin{ex}\normalfont
			The monomial basis of the $k$-algebra $R:=k[x, y]$ takes the form of Figure \ref{ fig corner points} when represented on the grid. If the monomial ideal $I$ of $k[x, y]$ is given by $I=\langle x^4, x^3y, x^2y^2, xy^3, y^5 \rangle$, then the quotient $k$-module $M=k[x, y]/I$ is $11$ dimensional and is generated by all elements in the Young diagram given in Figure \ref{ fig corner points}. The outside corner elements, which are circled red, generate $\mathfrak{R}(M)$, i.e., $\mathfrak{R}(M)=\langle x^3, x^2y, xy^2, y^4\rangle_k$ mod $I$ and this is the socle of $M$.
			
			\begin{figure}[H]
				
				$$	\begin{tikzpicture}[scale=0.75]
					
					\draw [-] (0, 2) -- (0, 8);
					\draw [-] (1, 2) -- (1, 8);
					\node(p) at (0.5, 2.5) {\vdots};
					\node(p) at (1.5, 4.5) {\vdots};
					\node(p) at (2, 4.3) {$\cdots$};
					\node(p) at (2.5, 5.5) {\vdots};
					\node(p) at (3, 5.3) {$\cdots$};
					\node(p) at (3.5, 6.5) {\vdots};
					\node(p) at (4, 6.3) {$\cdots$};
					\node(p) at (4.5, 7.5) {$\cdots$};

					\draw [-] (2, 5) -- (2, 8);
					\draw [-] (3, 6) -- (3, 8);
					\draw [-] (4, 7) -- (4, 8);
					\draw [-] (0, 3) -- (1, 3);
					
					\draw [-] (0, 4) -- (1, 4);
					\draw [-] (0, 5) -- (2, 5);
					
					\draw [-] (0, 6) -- (3, 6);
					\draw [-] (0, 7) -- (5, 7);
					\draw [-] (0, 8) -- (5, 8);
					
					\node (p) at (0.5, 7.5) {1};
					\node (p) at (1.5, 7.5) {$x$};
					\node (p) at (2.5, 7.5) {$x^2$};
					\node[circle, draw=red!80, inner sep=0pt, minimum size=28pt] ($x^3$) at (3.5, 7.5) {$x^3$};
					
					\node (p) at (0.5, 6.5) {$y$};
					\node (p) at (1.5, 6.5) {$xy$};
					\node[circle, draw=red!80, inner sep=0pt, minimum size=28pt] ($x^2y$) at (2.5, 6.5) {$x^2y$};

					\node (p) at (0.5, 5.5) {$y^2$};
					\node[circle, draw=red!80, inner sep=0pt, minimum size=28pt] ($xy^2$) at (1.5, 5.5) {$xy^2$};

					\node (p) at (0.5, 4.5) {$y^3$};
					
					\node[circle, draw=red!80, inner sep=0pt, minimum size=28pt] ($y^4$) at (0.5, 3.5) {$y^4$};
					
				\end{tikzpicture}$$
				\caption{Generators of $\mathfrak{R}(M)$ on a Young diagram.}\label{ fig corner points}.
			\end{figure}
		\end{ex}
		\begin{ex}\label{Rikard}
			\normalfont 	In general, $\mathfrak{R}(M)$ for $M\notin \mathfrak{C}$ need not be a submodule of $M$.
			Consider the $\mathbb{Z}$-module $M:=\mathbb{Z} \oplus \mathbb{Z}/p^2 \mathbb{Z}$, where $p$ is a prime number. The elements $(1, \bar{1})$ and $(1, \bar{0})$ of $M$  are torsion-free and therefore belong to $\mathfrak{R}(M)$. However, the element $(0, \bar{1})=(1, \bar{1}) - (1, \bar{0})$ of $M$ does not belong to $\mathfrak{R}(M)$ since $p^2(0, \bar{1})=(0, \bar{0})$ but $p(0, \bar{1}) \neq (0, \bar{0})$. This shows that in this case $\mathfrak{R}(M)$ is not a submodule of $M$.
		\end{ex}
		\begin{rem}\normalfont
			The elements in $S$ as given in Theorem \ref{thm soc(M)=R(M)} were described in \cite{wolff2016survival}, as the truly isolated monomials of a survival complex of a semigroup.
		\end{rem}
		\begin{paragraph}\noi
			An $R$-module $M$ is {\it coreduced} if for all $a\in R$, $aM=a^2M$.
		\end{paragraph}
		\begin{cor}\label{cor R(M) cored}\normalfont
			For any $M\in \mathfrak{C}$, $\mathfrak{R}(M)$ is coreduced.
		\end{cor}
		\begin{prf}
			For all $a\in \langle x_1, \cdots, x_n\rangle$, $a\mathfrak{R}(M)=a^2\mathfrak{R}(M)=0$, since $\mathfrak{R}(M)=\text{Soc}(M)=(0:_M\langle x_1, \cdots, x_n\rangle)$. If $a\in R\setminus \langle x_1, \cdots, x_n \rangle$, then $a=a_0+f(x_1, \cdots, x_n)$ where $a_0\in k$. Now, $a\mathfrak{R}(M)=a_0\mathfrak{R}(M)=\mathfrak{R}(M)$ and $a^2\mathfrak{R}(M)$\\$=a_0^2\mathfrak{R}(M)=\mathfrak{R}(M)$.
		\end{prf}
		\section{Reduced modules via inverse systems}
		\begin{paragraph}\noi
			The Macaulay inverse system is a powerful method for solving problems about Artinian local algebras of the form $R/I$. It was for instance used in problems such as weak Lefschetz property \cite{harbourne2011inverse}, in Waring's problem \cite{geramita1996waring}, and in the classiﬁcation of Artinian Gorenstein rings, \cite{Elias2012isomorphism}.
			In this section, we establish Macaulay inverse correspondences between reduced modules.
		\end{paragraph}
		\begin{paragraph}\noi
			Let $k$ be an algebraically closed field of characteristic zero.
			Let $R=k[x_1, \cdots, x_n]$, $V$ be the $k$-vector space $\langle x_1, \cdots, x_n \rangle$ and $P=\bigoplus_{i\geq 0}\text {Sym}^iV$, the standard graded polynomial in $n$-indeterminates over $k$. If $V^*$ is the $k$-vector space dual to $V$, then we have $V^*=\langle X_1, \cdots, X_n\rangle$ and $\Gamma=D^k(V^*)=\bigoplus_{i\geq 0} \text {Hom}_k(P_i, k)$ the graded $P$-module of graded $k$-linear homomorphisms from $P$ to $k$.
			It is well-known that; 
			$\Gamma \cong k[X_1, \cdots. X_n]$ the power divided ring and
			$\Gamma$ is an $R$-module via the following action which is called {\it apolarity action.} 
			$$R\circ \Gamma \longrightarrow \Gamma$$
		\end{paragraph}
		\begin{equation*}
			(x^\alpha, X^\beta)\mapsto x^\alpha \circ X^\beta:=
			\begin{cases}
				\frac{\beta!}{(\beta - \alpha)!}X^{\beta - \alpha} & \text{if}\quad \beta_i \geq \alpha_i, \text{for~all}~ i=1, \cdots, n\\
				0 & \text{otherwise}.
			\end{cases}
		\end{equation*}
		where $x^\alpha=x_1^{\alpha_1}, \cdots, x_n^{\alpha_n}$ and
		$X^\beta=X_1^{\beta_1}, \cdots, X_n^{\beta_n}$.
		\begin{paragraph}\noi
			Macaulay's correspondence is a special case of Matlis duality, which gives a one-to-one correspondence between ideals of $R$ and finitely generated submodules of $\Gamma$. For any ideal $I\subseteq R$, the dual of $I$ $$I^\perp:=\{m\in\Gamma~ |~ I\circ m=0\}$$ is a finitely generated submodule of $\Gamma$ called the {\it Macaulay's inverse system of $I$}. Conversely, if $W$ is a finitely generated submodule of $\Gamma$, then $$W^\perp=\{r\in R~|~r\circ W=0\}$$ is an ideal of $R$. So, to each Artinian local algebra $R/I$, we associate a finitely generated submodule $I^\perp$ of $\Gamma$. Conversely, if $W$ is a finitely generated submodule of $\Gamma$, then $R/W^\perp$ is a local Artinian algebra. We write $(R/I)^\vee=I^\perp$ and $W^\vee=R/W^\perp$ respectively. For more details about inverse systems, see for instance \cite{elias2021constructive,Elias2012isomorphism,emsalem1995inverse,geramita1996waring, harbourne2011inverse,cuong2018commutative, macaulay1994algebraic}.
			Just like in Section $2$, we consider the Artinian local algebra $R/I$ as an $R$-module $M\in \mathfrak{C}$.
		\end{paragraph}
		\begin{defn}\normalfont
			Let $I$ be an ideal of $R$. An $R$-module $M$ is {\it $I$-reduced} if for all $m\in M$, $I^2m=0$ implies that $Im=0$. Note that $M$ is reduced if it is $I$-reduced for all ideals $I$ of $R$.
		\end{defn}	
		\begin{paragraph}\noi
			We denote the largest $I$-reduced submodule of $M$ by $\mathfrak{R}_I(M)$.
		\end{paragraph} 
		\begin{lm}\label{Lm intersection of redu}\normalfont
			Let $R$ be a ring. For any $R$-module $M$, ~$\mathfrak{R}(M)=\bigcap \limits_{I\subseteq R} \mathfrak{R}_I(M)$.
		\end{lm}
		\begin{prf}
			For any ideal $I$ of $R$, every reduced $R$-module $M$ is $I$-reduced. Hence, $\mathfrak{R}(M) \subseteq \mathfrak{R}_I(M)$ for any ideal $I$ of $R$. So, $\mathfrak{R}(M)\subseteq \bigcap \limits_{I\subseteq R} \mathfrak{R}_I(M)$. Conversely, let $m\in \bigcap \limits_{I\subseteq R} \mathfrak{R}_I(M)$. Then for all ideals $I$ of $R$, $I^2m=0$ implies that $Im=0$.
			Thus $m\in \mathfrak{R}(M)$ and  $\bigcap \limits_{I\subseteq R}\mathfrak{R}_I(M)\subseteq \mathfrak{R}(M)$.
		\end{prf}
		
		\begin{thm}\normalfont\label{thm R(M)=k}
			Let $R=k[x_1, \cdots, x_n]$ and $\Gamma = k[X_1, \cdots. X_n]$. $\Gamma$ is an $R$-module via apolarity action. Define $\Gamma_i:= k[X_1, \cdots, X_i]$ a submodule of $\Gamma$ and $J_i:=\langle x_{i+1}, x_{i+2}, \cdots, x_n \rangle$ an ideal of $R$.
			\begin{itemize}
				\item [1.] If $J$ is an ideal of $R$ contained in $J_i$, then $\Gamma_i$ is $J$-reduced and $\Gamma_i \subseteq J^\perp$.
				\item [2.] If $W$ is a submodule of $\Gamma$ contained in $\Gamma_i$, then $W$ is $J_i$-reduced and $J_i \subseteq W^\perp$.
				\item [3.] $\mathfrak{R}_{J_i}(\Gamma)=\Gamma_i\subseteq J_i^\perp$.
				\item [4.] $\mathfrak{R}(\Gamma)=k$.
				
			\end{itemize}
		\end{thm}
		\begin{prf}
			\begin{itemize}
				\item [1.] Let $a\in J\subseteq J_i=\langle x_{i+1}, \cdots, x_n\rangle$. For any $a\in J_i,  a\circ \Gamma_i=0$. This is because $a$ consists of indeterminates $x_t$ such that $t>i$ for all indeterminates $X_i$ of $\Gamma_i$. In this case, $x_t\circ X_i=0$. It is also true that $x_t\circ k=0$, for all $t\geq 1$. So, $J\circ \Gamma_i=0$. Therefore $\Gamma_i$ is $J$-reduced and $\Gamma_i\subseteq (0:_\Gamma J)=J^\perp$.
				\item [2.] If $W$ is a submodule of $\Gamma_i\subseteq \Gamma$ and $J_i=\langle x_{i+1}, \cdots, x_n \rangle$. $W$ has the form $\sum\limits_j kX_j^j$ for some $j$ between $0$ and $i$. Just like in $1$ above $J_i\circ W=0$. So $W$ is $J_i$-reduced and $J_i\subseteq (0:_RW)=W^\perp$.
				\item [3.] Let $m$ be a monomial in the $R$-module $\Gamma$ and $a$ be a monomial in the ideal $J_i$ of $R$. If $a\circ m=0$, then either $a=x_t^p$ and $m=X_i^q$ such that $t>i$ or $a$ involves a term $x_j^t$ and $m$ involves a term $X_j^s$ such that $t>s$. However, $m\in \mathfrak{R}_{J_i}(\Gamma)$ if and only if $m$ is of the former type. To be of the former type is equivalent to having $m\in \Gamma_i$. $m$ in the latter case cannot be in $\mathfrak{R}_{J_i}(\Gamma)$. For if $t>s\neq 0$, then $x_j^t \circ X_j^s=0$ but $x_j\circ X_j^s\neq 0$. So, $\mathfrak{R}_{J_i}(\Gamma)=\Gamma_i$. Since for all monomials $m\in \Gamma_i$ we have some $a\in J_i$ such that $a\circ m=0$, $\Gamma_i\subseteq(0:_\Gamma J_i)=J_i^\perp$.
				\item [4.] Every nonzero $I$-reduced submodule of $\Gamma$ contains $k$. So, $$k\subseteq \bigcap\limits_{I\subseteq R}\{\Gamma_i~|~\Gamma_i~\text{is~an} ~I\text{-reduced~ submodule ~of}~ \Gamma\}.$$ Suppose $T:=\bigcap \limits_{I\subseteq R}\{\Gamma_i~|~\Gamma_i~\text{is~an} ~I\text{-reduced~ submodule ~of}~ \Gamma\}$ and $k$ is strictly contained in $T$. Then, by Lemma \ref{Lm intersection of redu} $T$ is a reduced submodule of $M$ and takes the form $k\bigoplus kX_1\bigoplus kX_2\bigoplus \cdots \bigoplus kX_n \bigoplus \cdots \bigoplus kX_1^{t_1}\bigoplus \cdots \bigoplus kX_n^{t_n}$, for some $t_i\in \mathbb{Z}^+$. Let $a=(x_1, \cdots, x_n)\in R$. $a^{t_i+1}\circ T=0$, but $a\circ T\neq 0$, contradicting the fact that $T$ is reduced. So, $k=T=\mathfrak{R}(\Gamma)$.

			\end{itemize}
		\end{prf}

		\begin{prop}\normalfont\label{prop R(M) dual}
			For any $M:=R/I\in \mathfrak{C}$ and a maximal ideal ${\bf m}$ of $R$, $$\mathfrak{R}(M)^\vee=\frac{I^\perp}{{\bf m}\circ I^\perp}.$$ 
		\end{prop}
		\begin{prf}
			Follows from Theorem \ref{thm soc(M)=R(M)} and \cite[Proposition 2.4.3]{cuong2018commutative}.\end{prf}
		
		\begin{ex}\label{dual reduced}\normalfont
			Let $R=k[x_1, x_2]$ and $I=\langle x_1^2, x_1x_2, x_2^3 \rangle $. If $M=R/I$, then $\mathfrak{R}(M)=\langle x_1, x_2^2\rangle_k $ mod $I$ and $I^\perp=k\bigoplus (kX_1 \bigoplus kX_2)\bigoplus kX_2^2$.
			$\mathfrak{R}(M)^\vee =\frac{I^\perp}{{\bf m}\circ I^\perp}=\langle X_1, X_2^2 \rangle_k$ mod $k\bigoplus kX_2$, which is a quotient of $I^\perp$. 
		\end{ex}
		\begin{defn}\normalfont
			Let $M=\bigoplus_{d\in \mathbb{N}}M_d$ be a graded $R$-module. We define the {\it Hilbert function} of $M$  $HF(M, -)~: \mathbb{N} \longrightarrow \mathbb{N}$ to be $HF(M, d):=\text{dim}_kM_d~ \text{for~ all}~ d\in \mathbb{N}$. Furthermore, we define the {\it Hilbert series} of $M$ as $$HS(M, t):=\sum_{d\in \mathbb{N}}HF(M, d)t^d.$$
		\end{defn}
		\begin{defn}\normalfont
			Let ${\bf m}$ be a maximal ideal of $R$. $X\in I^\perp$ is called an {\it outside corner} element if $x_i\circ X\in {\bf m}\circ I^\perp$ for all $1\leq i\leq n$. An element $X\in I^\perp$ is {\it inner} if it is not an outside corner element.
		\end{defn}
		
		\begin{thm}\normalfont\label{Thm duality}
			For any $M\in \mathfrak{C}$,
			\begin{itemize}
				\item [1.]	 a submodule $\mathfrak{R}(M)$ of $M$ is generated by monomials $x_1^{a_1} \cdots x_k^{a_k}$ mod $I$, $0\leq k\leq n$ if and only if $\mathfrak{R}(M)^\vee$, the quotient of $I^\perp$ is generated by the monomials $X_1^{a_1} \cdots X_k^{a_k}$ mod ${\bf m}\circ I^\perp$, where $a_1, \cdots, a_k$ are nonnegative integers;
				\item [2.] $HS(\mathfrak{R}(M), t)=HS(\mathfrak{R}(M)^\vee, t)$;
				\item [3.] $\mathfrak{R}(M)^\vee$ is the largest reduced quotient of $I^\perp$.
			\end{itemize}
		\end{thm}

		\begin{prf}
			\begin{itemize}
				\item [1.]
				It is well known that there is a one-to-one correspondence between the $R$-module $M=R/I$ and $I^\perp$, the Macaulay inverse system of $I$, \cite[IV]{macaulay1994algebraic}. It is also known that $HS(R/I, t)=HS(I^\perp, t)$, \cite[Proposition 2.3.3, Proposition 2.2.19]{cuong2018commutative}. Combinatorially, the generators of $R/I$ and $I^\perp$ can be represented in a Young diagram (YD). To distinguish between the two Young diagrams, we name the one for the former $\text{YD}_1$ and for the latter $\text{YD}_2$. It is easy to check that ${\bf m}\circ I^\perp$ is generated by all the inner elements of $\text{YD}_2$. It follows that $\mathfrak{R}(M)^\vee$, the quotient of $I^\perp$ by ${\bf m}\circ I^\perp$ is generated by all the elements at outside corners of the $\text{YD}_2$. However, this is in a one-to-one correspondence with elements at the outside corners of $\text{YD}_1$, which generate $\mathfrak{R}(M)$, see Theorem \ref{thm soc(M)=R(M)}. So, the generators of $\mathfrak{R}(M)$ are in one-to-one correspondence with the generators of $\mathfrak{R}(M)^\vee$.
				\item [2.] It is clear from $1$.
				\item [3.] Since $\mathfrak{R}(M)^\vee=I^\perp$ mod ${\bf m}\circ I^\perp=\langle \text{outside ~corner~elements~of}~\text{YD}_2 \rangle, x_i\circ \mathfrak{R}(M)^\vee \in {\bf m}\circ I^\perp$, whenever $x_i^2\circ \mathfrak{R}(M)^\vee \in {\bf m}\circ I^\perp$ for all integers $1\leq i \leq n$. The reduced quotients of $I^\perp$ are generated by elements at the outside corners of YD$_2$. Since $\mathfrak{R}(M)^\vee$ is generated by all the elements at the outside corners of YD$_2$, it is the largest reduced quotient of $I^\perp$.
			\end{itemize}
		\end{prf}
		
		\begin{ex}\normalfont
			Consider the $R$-module $M=\frac{k[x, y]}{\langle x^4, x^3y, y^2 \rangle }$. $\mathfrak{R}(M)=\langle x^3, x^2y\rangle_k$ mod $\langle x^4, x^3y, y^2 \rangle$ and  $I^\perp=k\bigoplus kX\bigoplus kY\bigoplus kX^2\bigoplus kXY \bigoplus kX^3 \bigoplus kX^2Y$. So, $\mathfrak{R}(M)^\vee\\=kX^3 \bigoplus kX^2Y$ mod $(k\bigoplus kX \bigoplus kY \bigoplus kX^2 \bigoplus kXY)$. The $\text{YD}_1$ and $\text{YD}_2$ associated to $M$ and $I^\perp$ are respectively given in Figure \ref{Fig YD_1 and YD_2}.
			
			\begin{figure}[H]
				$$	\begin{tikzpicture} [scale=1]
					\fill [yellow ]  (2, 6) --  (5, 6) rectangle (2, 5) -- (5, 5);	\fill [yellow ]  (2, 5) --  (4, 5) rectangle (2, 4) -- (4, 4);
					\fill [orange ]  (5, 6) --  (6, 6) rectangle (5, 5) -- (6, 5);	\fill [orange ]  (4, 5) --  (5, 5) rectangle (4, 4) -- (5, 4);
					
					\draw[-](2, 4) -- (5, 4);
					\draw[-](2, 5) -- (6, 5);
					\draw[-] (2, 4) -- (2, 5);
					\draw[-] (5, 4) -- (5, 5);
					\draw[-] (3, 4) -- (3, 5);
					\draw[-] (4, 4) -- (4, 5);
					\draw[-] (2, 5) -- (2, 6);
					\draw[-] (3, 5) -- (3, 6);
					\draw[-] (4, 5) -- (4, 6);
					\draw[-] (5, 5) -- (5, 6);
					\draw[-] (2, 6) -- (6, 6);
					\draw[-] (6, 5) -- (6, 6);

					\node (P) at (2.5, 5.5) {1 };
					\node (P) at (3.5, 5.5) {$x$ };
					\node (p) at (4.5,5.5) {$x^2$};
					\node (p) at (5.5,5.5) {$x^3$};
					
					\node (P) at (2.5, 4.5) {$y$ };
					\node (p) at (3.5,4.5) {$xy$};
					\node (p) at (4.5,4.5) {$x^2y$};

					\fill [yellow ]  (8, 6) --  (11, 6) rectangle (8, 5) -- (11, 5);	\fill [yellow ]  (8, 5) --  (10, 5) rectangle (8, 4) -- (10, 4);
					
					\fill [orange ]  (11, 6) --  (12, 6) rectangle (11, 5) -- (12, 5);	\fill [orange ]  (10, 4) --  (11, 4) rectangle (10, 5) -- (11, 5);
					
					\draw[-](8, 6) -- (12, 6);
					\draw[-](8, 5) -- (12, 5);
					\draw[-] (8, 4) -- (11, 4);
					\draw[-] (8, 4) -- (8, 6);
					\draw[-] (9, 4) -- (9, 6);
					\draw[-] (10, 4) -- (10, 6);
					\draw[-] (11, 4) -- (11, 6);
					\draw[-] (12, 5) -- (12, 6);

					\node (P) at (8.5, 4.5) {$Y$ };
					\node (P) at (8.5, 5.5) {$1$ };
					\node (p) at (9.5,5.5) {$X$};
					\node (p) at (10.5,5.5) {$X^2$};
					\node (p) at (11.5,5.5) {$X^3$};
					
					\node (P) at (9.5, 4.5) {$XY$ };
					\node (p) at (10.5,4.5) {$X^2Y$};	
				\end{tikzpicture}$$
				\caption{$\text{YD}_1$ and $\text{YD}_2$ of $M$ and $I^\perp$ respectively.}\label{Fig YD_1 and YD_2}
			\end{figure}
		\end{ex}
		Define
		$\text{\Large $\mathfrak{D}$}:=\bigg\{ I^\perp\leq \Gamma: \text{for~ an~ ideal}~ I ~\text{of}~ R~\text{such~that}~R/I\in \mathfrak{C}\bigg \}$.
		
		\begin{cor}\normalfont\label{cor properties of perp and vee}
			The following statements hold for any $M\in \mathfrak{C}$ and $I^\perp \in \mathfrak{D}$.
			\begin{itemize}
				\item [1.] $\mathfrak{R}(M)^\vee$ is both a reduced and a coreduced quotient of $I^\perp.$
				\item [2.] $\mathfrak{R}(M^\vee)=\mathfrak{R}(I^\perp)=k$ (the only reduced submodule of $\Gamma$).
				\item [3.] $\mathfrak{R}(I^\perp)^\vee\cong k$ (the only reduced quotient of $M$).
				\item [4.] $\mathfrak{R}((I^\perp)^\vee)=\mathfrak{R}(M)$ (is both a reduced and a coreduced submodule of $M$).
			\end{itemize}
		\end{cor}
		\begin{prf}
			\begin{itemize}
				\item [1.] Since $\mathfrak{R}(M)^\vee=I^\perp$ mod ${\bf m}\circ I^\perp=\langle \text{outside ~corner~elements~of}~\text{YD}_2 \rangle, x_i\circ \mathfrak{R}(M)^\vee \in {\bf m}\circ I^\perp$ and $x_i^2\circ \mathfrak{R}(M)^\vee \in {\bf m}\circ I^\perp$, for all $i\in \mathbb{Z}^+$ which implies $x_i\circ \mathfrak{R}(M)^\vee=x_i^2\circ \mathfrak{R}(M)^\vee$ for all $i\in \mathbb{Z}^+$. This shows that $\mathfrak{R}(M)^\vee$ is a coreduced $R$-module. For the reduced part see, Theorem \ref{Thm duality} (3).
				\item [2.] $\mathfrak{R}(M^\vee)=\mathfrak{R}(I^\perp)=I^\perp\cap\mathfrak{R}(\Gamma)=I^\perp \cap k=k$.
				\item [3.] $\mathfrak{R}(I^\perp)^\vee=k^\vee= R/k^\perp =R/{\bf m}\cong k$.
				\item [4.] Follows by definition.
			\end{itemize}
		\end{prf}
		\begin{prop}\normalfont
			Let $R=k[x_1, \cdots, x_n]$ and $I$ be an ideal of $R$ which is homogeneous of degree $n+1$. If $M:=R/I$, then 
			\begin{itemize}
				\item [1.] $\mathfrak{R}(M)=(x_1, \cdots, x_n)^nM=M_n$ (the homogeneous part of $M$ with the highest degree).
				\item [2.] $\mathfrak{R}(\Gamma)=(x_1, \cdots, x_n)^n\circ I^\perp=k$ (the homogeneous part of $I^\perp $ with the least degree).
			\end{itemize}
		\end{prop}
		\begin{prf}
			\begin{itemize}
				\item [1.] Successive multiplication of an element $(x_1, \cdots, x_n)$ with $M:=R/I$ a graded $R$-module eliminates on each multiplication the homogeneous part of $M$ with the lowest degree. Therefore multiplying $n$-times leaves only the homogeneous part of $R/I$ with the highest degree. This is however, the submodule of $M$ generated by the outside corner elements of $M$. So, $\mathfrak{R}(M)=(x_1, \cdots, x_n)^nM$.
				\item [2.] $I^\perp$ is a finite dimensional graded $R$-module. Successive multiplication of $I^\perp$ by the element $(x_1, \cdots, x_n)$ via apolarity action eliminates the homogeneous part with the highest degree. Multiplying $n$ times leaves only the homogeneous part of $I^\perp$ with degree $0$, which is $k$. By Theorem \ref{thm R(M)=k}, $\mathfrak{R}(\Gamma)=k$. So, $\mathfrak{R}(\Gamma)=(x_1, \cdots, x_n)^n\circ I^\perp=k$.
			\end{itemize}			
		\end{prf}
		
		\begin{cor}\normalfont
			For any $I^\perp\in \mathfrak{D}$,  $\mathfrak{R}(I^\perp)=\text{Soc}(I^\perp)$.
		\end{cor}
		\begin{prf}
			For any monomial ideal $I$ of $R$, the corresponding submodule $I^\perp$ of $\Gamma$ contains $k$ and this is the only submodule of $\Gamma$ which is simple, hence $\text{Soc}(I^\perp)=k$. Moreover, by Corollary \ref{cor properties of perp and vee} $k$ is the largest reduced submodule of $I^\perp$. Thus, $\mathfrak{R}(I^\perp)=\text{Soc}(I^\perp)$.
		\end{prf}
		\begin{paragraph}\noi
			Let $\mathfrak{D}_\text{red}$ be the full subcategory of $R$-Mod
			consisting of all reduced $R$-modules $N$ defined under apolarity action such that there exists a surjection map $I^\perp \twoheadrightarrow N$ for some $I^\perp \in \mathfrak{D}$. Then $\mathfrak{D}_\text{red}$ is dual to $ \mathfrak{C}_\text{red}$	and coincides with all semisimple modules $T$ defined under apolarity action for which the surjection $I^\perp \twoheadrightarrow T$ exists for some $I^\perp \in \mathfrak{D}$.	
		\end{paragraph}
		\section{The symmetries exhibited}
		\begin{paragraph}\noi
			This section is mostly complementary to papers \cite{ssevviiri2022applications,ssevviiri2023applications}. In particular, it serves three purposes:
			1) it demonstrates that $\mathfrak{C}_\text{red}$, provides a concrete and accessible example of modules: a) in the TTF class $T_I$ constructed in \cite[Theorem 3.1]{ssevviiri2023applications} for all ideals $I$ of $R$, b) for which the Matlis-Greenlees-May  Equivalence proved in \cite[Theorem 4.3]{ssevviiri2022applications} and later showed to be an equality in \cite[Proposition 4.2]{ssevviiri2023applications} holds for all ideals $I$ of $R$; 2) it exhibits:
			a) commutativity between taking torsion theories and taking Matlis duality, b) symmetries between the $R$-modules $M, \Gamma_I(M), \Lambda_I(M), M/IM$, and $(0:_MI)$ and their associated Matlis duals; 3) gives new results and a summary of  several results studied in \cite{ssevviiri2022applications, ssevviiri2023applications}.
		\end{paragraph}
		\begin{paragraph}\noi
			For an $R$-module $M$, define $\Gamma_I(M):=\underset{k}\varinjlim \text{Hom}(R/I^k, M)$ and \\$\Lambda_I(M):=\underset{k}\varprojlim M/I^kM$.
		\end{paragraph}
		\begin{paragraph}\noi
			Figure \ref{Fig table} summarizes the Macaulay inverse correspondences of reduced modules studied in Section 3.
		\end{paragraph}
		\begin{figure}[H]
			$$	\begin{tikzpicture}[scale=0.6]
				\draw[-](0,16) -- (25,16);
				\draw[-](0,14) -- (25,14);
				\draw[-](0,12) -- (25,12);
				\draw[-](0,10) -- (25,10);
				\draw[-](0,8) -- (25,8);
				\draw[-](0,6) -- (25,6);
				\draw[-](0,4) -- (25,4);
				\draw[-](0,1.8) -- (25,1.8);
				\draw[-](0,0) -- (25,0);
				\draw[-](0,-2) -- (25,-2);
				\draw[-](0,-4.5) -- (25,-4.5);
				
				\draw[-](0,16) -- (0,-4.5);
				\draw[-](25,16) -- (25,-4.5);

				\node (P) at (2.5, 15.5) {Modules under};
				\node (P) at (2.5, 14.5) { usual action};
				\node (P) at (9.7, 15.5) {Modules under};
				\node (P) at (9.7, 14.5) { apolarity action};
				\node (P) at (17, 15) {Remarks};
				
				\node (P) at (0.5, 13) {1.};
				
				\node (P) at (3, 13) {$M=R/I\in \mathfrak{C}$};
				\node (P) at (9.3, 13) {$I^\perp \in\mathfrak{D}$};
				\node (P) at (7, 13) {$\longleftrightarrow$ };
				
				\node (P) at (18, 13) { the two have the same dimension};

				\node (P) at (0.5, 11) {2.};
				
				\node (P) at (3, 11) {$\mathfrak{R}(M)$ };
				\node (P) at (7, 11) {$\longleftrightarrow$ };
				\node (P) at (9, 11) {$\frac{I^\perp}{{\bf m}\circ I}$ };
				\node (P) at (17.8, 11.5) { both are reduced, coreduced and};
				\node (P) at (18.5, 10.5) {generated by outside corner elements};
				\node (P) at (0.5, 9) {3.};
				
				\node (P) at (3.5, 9) {$\frac{k[x_1,\cdots, x_n]}{\langle x_1,\cdots, x_1\rangle}\cong k\in \mathfrak{C}$ };
				\node (P) at (7, 9) {$\longleftrightarrow$ };
				\node (P) at (9.3, 9) { $k\in \mathfrak{D}$};	
				\node (P) at (18.6, 9.5) { the only reduced modules in $\mathfrak{C}$ and $\mathfrak{D}$};
				\node (P) at (14.7, 8.5) { respectively};
				
				\node (P) at (0.5, 7) {4.};
				\node (P) at (3, 7) {$M/\mathfrak{R}(M)$ };
				\node (P) at (7, 7) {$\longleftrightarrow$ };
				\node (P) at (9.3, 7) { ${\bf m}\circ I^\perp$};	
				\node (P) at (17.2, 7) { generated by inner elements};
				
				\node (P) at (0.5, 5) {5.};
				\node (P) at (3, 5) {$\mathfrak{R}(M) \hookrightarrow  M$ };
				\node (P) at (7, 5) {$\longleftrightarrow$ };
				\node (P) at (10, 5) {$I^\perp  \twoheadrightarrow  \frac{I^\perp}{{\bf m}\circ I^\perp}$};	
				\node (P) at (17.7, 5.5) { an embedding and a surjection};
				\node (P) at (14.9, 4.5) { respectively};	
				
				\node (P) at (0.5, 3) {6.};
				\node (P) at (3, 3.5) {$M\twoheadrightarrow M/\bar{{\bf m}},$} ;
				\node (P) at (3.2, 2.4) { $\bar{{\bf m}}=\frac{\langle x_1,\cdots, x_1\rangle}{I}$};
				
				\node (P) at (7, 3) {$\longleftrightarrow$ };
				\node (P) at (10, 3) {$\mathfrak{R}(I^\perp) \hookrightarrow I^\perp$};	
				\node (P) at (17.8, 3.5) { a surjection and an embedding};
				\node (P) at (14.9, 2.5) { respectively};
				
				\node (P) at (0.5, 1) {7.};
				\node (P) at (3, 1) {$\mathfrak{C}_\text{red}$ };
				\node (P) at (7, 1) {$\longleftrightarrow$ };
				\node (P) at (10, 1) {$\mathfrak{D}_\text{red}$ };	
				\node (P) at (18.5, 1) { both consist of semisimple modules};
				
				\node (P) at (0.5, -1) {8.};
				\node (P) at (3.5, -1) {$\mathfrak{R}(M)=\text{Soc}(M)$ };
				\node (P) at (7, -1) {$\longleftrightarrow$ };
				\node (P) at (10.6, -1) {$\mathfrak{R}(I^\perp)=\text{Soc}(I^\perp)$ };	
				\node (P) at (19, -0.5) { in both cases the  reduced submodule};
				\node (P) at (17, -1.5) { and the socle coincide};
				
				\node (P) at (0.5, -3) {9.};
				\node (P) at (3.8, -3) {$\text{Soc}(M)=(0:_M{\bf m})$ };
				\node (P) at (7.5, -3) {$\longleftrightarrow$ };
				\node (P) at (11, -3) {$\text{Soc}(I^\perp)=(0:_\Gamma{\bf m})$ };	
				\node (P) at (19.5, -3) { Socle is the annihilating submodule};
				\node (P) at (19, -4) {  by ${\bf m}$ of $M$ and $\Gamma$ respectively};

			\end{tikzpicture}$$
			\caption{The summary of Macaulay inverse correspondences about reduced modules.}\label{Fig table}
			
		\end{figure}

		\begin{defn}\normalfont
			An $R$-module $M$ is {\it $I$-torsion} (resp. {\it $I$-complete, $I$-reduced} and {\it $I$-coreduced} ) if $\Gamma_I(M)\cong M$ (resp. $\Lambda_I(M)\cong M, \Gamma_I(M)\cong \text{Hom}(R/I, M)$ and $\Lambda_I(M)\cong R/I\otimes M)$.
		\end{defn}
		\begin{defn}\normalfont
			A {\it torsion theory} of an abelian category $\mathcal{C}$, is a pair $(T, F)$ of full subcategories of $\mathcal{C}$ such that 
			$\text{Hom}(T, F)=0$ and for all $M\in \mathcal{C}$, there exists a short exact sequence $$0\rightarrow M_T\rightarrow M\rightarrow M_F\rightarrow 0$$ with $M_T\in T$ and $M_F\in F$. In this case, we call $T$ a {\it torsion class} and $F$ a {\it torsionfree class}. A class $\mathcal{L}$ of an abelian category $\mathcal{C}$ is a {\it torsion-torsionfree} (TTF) class  if it is both a torsion and a torsion-free class.
		\end{defn}
		\begin{paragraph}\noi
			Let $\mathcal{A}_I$ (resp. $\mathcal{B}_I$) be an abelian full subcategory of $R$-Mod consisting of $I$-reduced (resp. $I$-coreduced) $R$-modules. In \cite[Theorem 3.1]{ssevviiri2023applications}; it was shown that the class
			$$T_I=\{M\in \mathcal{A}_I~|~\Gamma_I(M)=M\}$$ is a TTF with the associated torsion-torsionfree triple $(\mathfrak{T}_I, T_I, \mathcal{F}_I)$, where $\mathcal{F}_I=\{M\in \mathcal{A}_I:\Gamma_I(M)=0\}$ and $\mathfrak{T}_I=\{M\in \mathcal{B}_I:IM=M\}$.
		\end{paragraph}
		\begin{ex}\normalfont
			$\mathfrak{C}_\text{red}\subseteq T_I$, for all ideals $I$ of $R$. Secondly, any $M=R/I\in \mathfrak{C}$ is also contained in $T_I$. In these two cases the modules are $I$-reduced and $I$-torsion if and only if they are $I$-coreduced and $I$-complete, \cite[Proposition $4.2$]{ssevviiri2023applications}.
		\end{ex}
		\begin{prop}\normalfont\label{thm TTF}
			Let $M$ be an $R$-module which is both $I$-reduced and $I$-coreduced. Consider the TTF triple $(\mathfrak{T}_I, T_I, \mathcal{F}_I)$  and the Matlis dual $\text{Hom}_R(-, E)$, then
			
			\begin{enumerate}
				\item $M\in T_I$ if and only if $\text{Hom}_R(M, E)\in T_I$;
				\item let in addition $I$ be finitely generated, $M\in \mathcal{F}_I$ if and only if $\text{Hom}_R(M, E)\in \mathfrak{T}_I$;
				\item $M\in \mathfrak{T}_I$ if and only if $\text{Hom}_R(M, E)\in \mathcal{F}_I$.
			\end{enumerate}
		\end{prop}
		\begin{prf}
			\begin{enumerate}
				\item  $\Gamma_I(M)\cong M \Leftrightarrow \Lambda_I(M)\cong M\Leftrightarrow\text{Hom}_R(\Lambda_I(M), E)\cong \text{Hom}_R(M, E)\Leftrightarrow \Gamma_I(\text{Hom}_R(M,E))\cong \text{Hom}_R(M,E).$ The first equivalence is due to Proposition \cite[Proposition 4.2]{ssevviiri2023applications}. The  second equivalence holds  because the functor $\text{Hom}_R(-, E)$ preserves and reflects isomorphisms. The third equivalence is due to \cite[Proposition 5.3]{ssevviiri2022applications}.
				
				\item $\Gamma_I(M)=0 \Leftrightarrow\text{Hom}_R(\Gamma_I(M), E)\cong 0 \Leftrightarrow \Lambda_I(\text{Hom}_R(M,E))\cong 0$. The part of the first equivalence is trivial; its reverse is due to the fact that $\text{Hom}_R(-, E)$ reflects zero since $E$ is an injective cogenerator.  The second equivalence is due to \cite[Proposition 5.5]{ssevviiri2022applications}.
				
				\item First note that $IM=M$ if and only if $ \Lambda_I(M)\cong 0\Leftrightarrow\text{Hom}_R(\Lambda_I(M), E)\cong 0 \Leftrightarrow \Gamma_I(\text{Hom}_R(M,E))\cong 0.$ The first equivalence  holds because $\text{Hom}_R(-, E)$ reflects and preserves zeros. The second equivalence is due to \cite[Proposition 5.3]{ssevviiri2022applications}.
			\end{enumerate}
			
		\end{prf}
		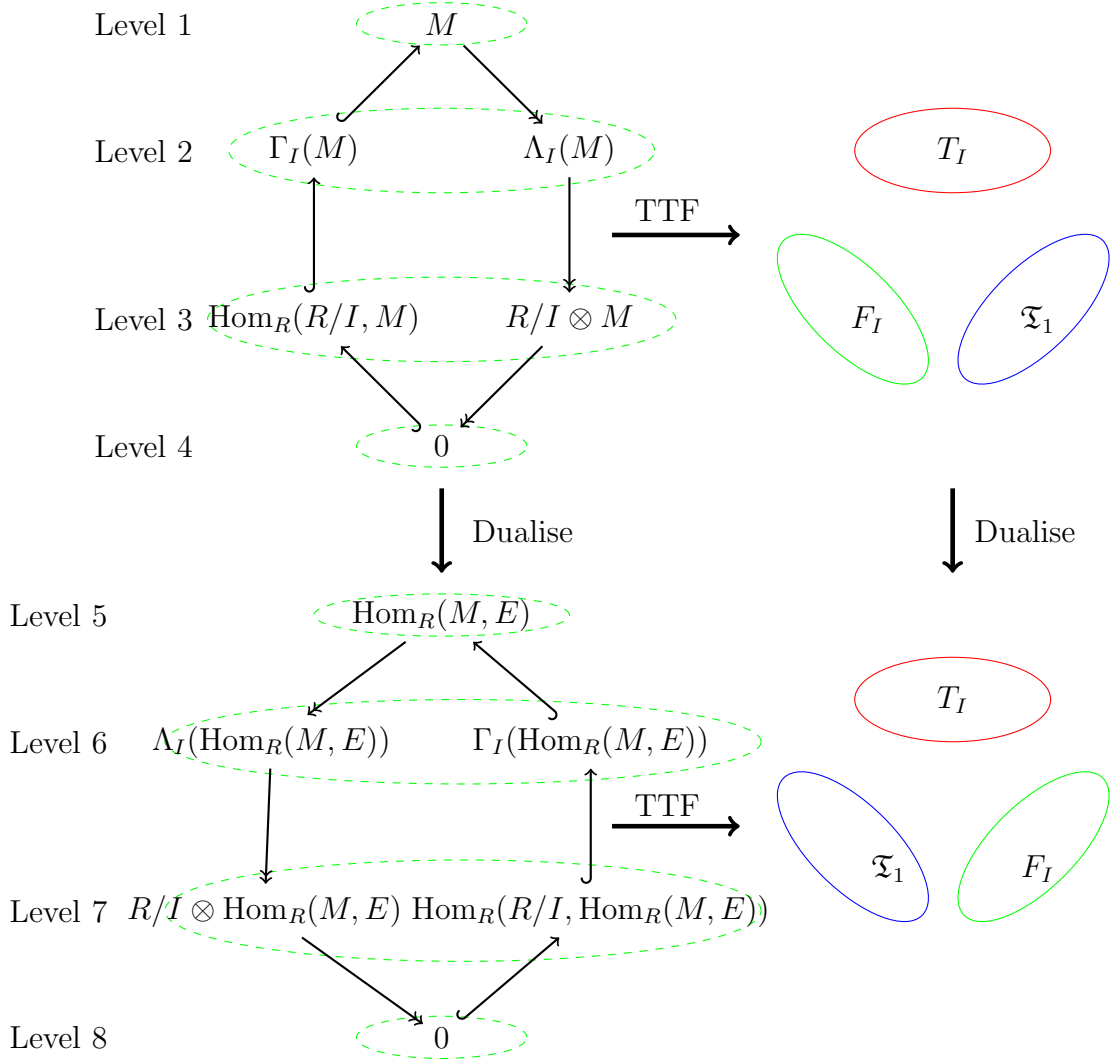
\begin{figure}[H]
			
			\begin{tikzpicture}[scale=0.56]
				\draw[green, dashed](-4,4)  ellipse (5.5cm and 1 cm);
				\draw[green, dashed](-4,8)  ellipse (5.0cm and 1 cm);
				\draw[green, dashed](-4,11)  ellipse (2cm and 0.5 cm);
				\draw[green, dashed](-4, 1)  ellipse (2cm and 0.5 cm);
				\node (A) at (-7, 8) {$\Gamma_I(M)$};
				
				\node (B) at (-7, 4) {$\text{Hom}_R(R/I, M)$};
				\node (C) at (-1, 4) {$R/I\otimes M$};
				\node (D) at (-1, 8) {$\Lambda_I(M)$};
				\node (E) at (-4, 1) { $0$};
				\node (F) at (-4, 11) { $M$};
				\draw[thick,  right hook ->] ( E) edge (B) (B) edge (A) (A) edge (F);
				
				\draw[thick, ->>] (F) edge (D) (D) edge (C)
				(C) edge (E);

				\draw[red] (8, 8) ellipse (2.3cm and 1.0 cm);
				\node at (8, 8) {$T_I$};

				\draw[rotate = 45, blue] (10, -4) ellipse (2.3cm and 1.0 cm);
				\node at (10, 4) {$\mathfrak{T}_1$};

				\draw[rotate = -45, green] (1, 7) ellipse (2.3cm and 1.0 cm);
				\node at (6, 4) {$F_I$};
				
				\draw[green, dashed](-3.5,-10)  ellipse (7.0cm and 1.2 cm);  
				\draw[green, dashed](-3.5,-6)  ellipse (7cm and 1 cm);
				\draw[green, dashed](-4,-3)  ellipse (3cm and 0.5 cm);
				\draw[green, dashed](-4,-13)  ellipse (2cm and 0.5 cm);
				\node (A) at (-8, -6) {$\Lambda_I(\text{Hom}_R(M, E))$};
				\node (B) at (-8.16, -10) {$R/I\otimes \text{Hom}_R(M, E)$};
				\node (C) at (-0.5, -10) {$\text{Hom}_R(R/I,  \text{Hom}_R(M, E))$};
				\node (D) at (-0.5, -6) {$\Gamma_I(\text{Hom}_R(M, E))$};
				\node (E) at (-4, -13) { $0$};
				\node (F) at (-4, -3) { $\text{Hom}_R(M, E)$};
				\draw[thick, ->>]
				(B) edge (E)
				(A) edge (B) (F) edge (A);
				\draw[thick, right hook->]
				(D) edge (F) (C) edge (D) (E) edge (C);
				
				\node (P) at (-13.0, -3) { Level 5};
				\node (Q) at (-13.0, -6) {Level 6};
				\node (R) at (-13.0, -10) { Level 7};
				\node (S) at (-13.0, -13) { Level 8};

				\draw[red] (8, -5) ellipse (2.3cm and 1.0 cm);
				\node at (8, -5) {$T_I$};

				\draw[rotate = 45, green] (1, -13) ellipse (2.3cm and 1.0 cm);
				\node at (10, -9) {$F_I$};

				\draw[rotate = -45, blue] (10, -2) ellipse (2.3cm and 1.0 cm);
				\node at (6.5, -9) {$\mathfrak{T}_1$};
				
				\draw[->, line width = 0.61mm] (0, 6) -- (3, 6);
				\draw[->, line width = 0.61mm] (0, -8) -- (3,-8);
				\draw[->, line width = 0.61mm] (8, 0) -- (8, -2);
				\draw[->, line width = 0.61mm] (-4, 0) -- (-4, -2);

				\node (P) at (-11.0, 11) { Level 1};
				\node (Q) at (-11.0, 8) {Level 2};
				\node (R) at (-11.0, 4) { Level 3};
				\node (S) at (-11.0, 1) { Level 4};
				
				\node (T) at (1.3, -7.5) {TTF};
				
				\node (U) at (1.3, 6.6) {TTF};
				
				\node (W) at (-2.1, -1) {Dualise};
				
				\node (X) at (9.7, -1) {Dualise};
				
			\end{tikzpicture}
			\caption{Symmetries summarised.}\label{fig summary}  
		\end{figure}

		\begin{paragraph}\noi
			For any $M\in \mathfrak{C}_\text{red}$, levels 1, 2 and 3  of Figure \ref{fig summary} collapse to one level, i.e., $\Gamma_I(M)\cong \text{Hom}_R(R/I, M)\cong M$ and $\Lambda_I(M)\cong R/I\otimes M \cong M$. In this case, $M$ belongs to the torsion class $T_I$ (in the red ellipse). By Matlis duality levels $5, 6$ and $7$ also collapse to one level and $\text{Hom}(M, E)$ also belongs to $T_I$.
		\end{paragraph}
		\begin{paragraph}\noi
			If $\Gamma_I(M)=0$, then the left hand side of levels $2, 3$ and $4$ collapse to one level. In this case $M\in \mathcal{F}_I$ (the green ellipse). However, on taking the Matlis dual of the module $\Gamma_I(M)$, $\text{Hom}(R/I, M)$ and $0$ we get the left hand side of levels $6, 7$ and $8$ collapse to one level in which case $\text{Hom}(M, E) \in \mathfrak{T}_I$ (the blue ellipse). 
		\end{paragraph}
		
		\begin{paragraph}\noi
			If $\Lambda_I(M)=0$, the right hand side of levels $2, 3$ and $4$ collapse to one level. So, $M\in \mathfrak{T}_I$ (the blue ellipse). Taking the Matlis dual of the modules $\Lambda_I(M), R/I\otimes M$ and $0$, one gets the right hand side of levels $6, 7$ and $8$ also collapse to just one level and $\text{Hom}(M, E)\in \mathcal{F}_I$ (the green ellipse).
		\end{paragraph}

		\begin{paragraph}\noi
			The following results can also be interpreted using Figure \ref{fig summary}.
			\begin{itemize}
				\item [1.] Let $I$ be finitely generated ideal of a ring $R$. If $M$ is an $I$-reduced and $I$-torsion $R$-module then $\text{Hom}_R(M, E)$ is $I$-complete, \cite[Corollary 5.6 (2)]{ssevviiri2022applications}.
				\item [2.] $M$ is $I$-coreduced, then $\text{Hom}_R(M, E)$ is $I$-reduced, \cite[Proposition 2.6 (1)]{ssevviiri2022applications}.
				\item[3.] Let $R$ be a Noetherian ring. $M$ is $I$-reduced if and only if $\text{Hom}_R(M, E)$ is $I$-coreduced, \cite[Proposition 5.1]{ssevviiri2022applications}.
				\item [4.] Let $I$ be any ideal of a ring $R$ and $N$ any $R$-module. If $M$ is $I$-coreduced and $I$-complete, then $\text{Hom}_R(M, N)$ is $I$-torsion, \cite[Corollary 5.4 (2)]{ssevviiri2022applications}.
				\item [5.] Let $R$ be a Noetherian ring, and $\mathcal{A}_I$ be an abelian full subcategory of $R$-Mod consisting of $I$-reduced $R$-modules. The $I$-torsion free modules in $\mathcal{A}_I$ coincide with the $I$-coreduced $R$-modules $M$ for which $IM=M$, i.e., $\mathfrak{T}_I=\mathcal{F}_I$ (the green and the blue ellipses coincide), \cite[Proposition 3.3 (2)]{ssevviiri2023applications}.
			\end{itemize}
		\end{paragraph}

		\begin{paragraph}\noi
			The two processes namely; (1) of taking a TTF and then dualising and
			(2) of dualising and then taking the TTF are commutative. The arrow between the ellipses summarises Proposition \ref{thm TTF}.
		\end{paragraph}

		\section{Modules in $\mathfrak{C}$ satisfy the radical formula}
		
		\begin{paragraph}\noi
			A submodule $P$ of an $R$-module $M$ is {\it prime} (resp. {\it semiprime}) if the $R$-module $M/P$ is prime (resp. reduced). The prime radical $\beta(M)$ (resp. semiprime radical, $\mathcal{S}(M)$) of $M$ is the intersection of all prime  (resp. semiprime) submodules of $M$. If $N$ is a submodule of $M$, then $\beta(N)$ (resp. $\mathcal{S}(N)$) denotes the intersection of all prime (resp. semiprime) submodules of $M$ containing $N$.
			For a submodule $N$ of an $R$-module $M$, we define the envelope $E_M(N)$ of $N$ as the set $$E_M(N):=\{rm~|~ r\in R, m\in M~\text{and}~r^km\in N~\text{for some}~ k\in \mathbb{Z}^+\}.$$ Denote the submodule of $M$ generated by $E_M(N)$ by $\langle E_M(N)\rangle$. A submodule $N$ of $M$ is said {\it to satisfy the radical formula} (or  s.t.r.f for short)  if $\langle E_M(N)\rangle=\beta(N)$. A module $M$ is said to {\it satisfy the radical formula}  if every submodule $N$ of $M$ s.t.r.f. Modules that s.t.r.f have been studied in  \cite{azizi2007radical, azizi2009radical, Jenkins1992,lee2019semiprime,leung1997,mccasland1991radicals,pusat1986,  sharif2002,david2015}  among others. Lastly, we define the Jacobson radical $\mathcal{J}(M)$ of $M$ as the intersection of all maximal submodules of $M$.
		\end{paragraph}
		\begin{paragraph}\noi
			$\mathcal{S}(M)$ (resp. $\mathcal{J}(M))$ is dual to $\mathfrak{R}(M)$ (resp. $\text{Soc}(M)$). As expected, (since $\mathfrak{R}(M)=\text{Soc}(M)~ \text{for}~M\in \mathfrak{C}$) the submodules $\mathcal{S}(M)$ and $\mathcal{J}(M)$ also coincide whenever $M\in \mathfrak{C}$.
		\end{paragraph}
		\begin{thm}\normalfont
			For any $M\in \mathfrak{C}$,
			\begin{itemize}
				\item [1.] $\langle E_M(0)\rangle=\mathcal{S}(M)=\beta(M)=\mathcal{J}(M)$,
				\item [2.] $M$ s.t.r.f.
			\end{itemize}
		\end{thm}
		\begin{prf}
			\begin{itemize}
				\item [1.]
				The inclusion $\langle E_M(0)\rangle\subseteq \mathcal{S}(M)\subseteq \beta(M)\subseteq \mathcal{J}(M)$ is well-known but also easy to see, since we have the following implications between submodules,\\ maximal $\Rightarrow $ prime $\Rightarrow$ semiprime. For $\langle E_M(0)\rangle \subseteq \mathcal{S}(M)$, see for instance \cite[section 3]{lee2019semiprime}. Since $M$ has only one maximal submodule; namely $\bar{{\bf m}}=\langle x_1, \cdots, x_n\rangle/I$, $\mathcal{J}(M)=\bar{{\bf m}}$.
				Note also that $\langle E_M(0)\rangle=\langle x_1, \cdots, x_n\rangle M=\bar{{\bf m}}$, so $\langle E_M(0)\rangle=\mathcal{J}(M)$.
				\item [2.]
				From part $1)$ $\langle E_M(0)\rangle=\beta(M)$, i.e., the zero submodule of $M$ s.t.r.f. By \cite[Theorem 1.5]{ mccasland1991radicals} given an epimorphism $f:M\rightarrow M/N$ for any submodule $N$ of $M$, the zero submodule s.t.r.f if and only if every submodule $N$ of $M$ s.t.r.f.
			\end{itemize}
		\end{prf}
		\begin{cor}\normalfont
			Any module $M\in \mathfrak{C}$, has exactly one semiprime submodule; namely, $\langle x_1, \cdots, x_n \rangle/I$.
		\end{cor}
		
		\subsection*{Acknowledgment}
		\begin{paragraph}\noi
			The authors acknowledge support from the International Science Program, EMS-Simons for Africa and the Eastern Africa Algebra Research Group. The authors are grateful to Rikard Bøgvad for pointing out Example \ref{Rikard} to them and to Dominic Bunnett, Alexandru Constantinescu and Balazs Szendroi for the comments. The third author is grateful to the fourth author for facilitating his visit to Makerere University.
		\end{paragraph}

		\begin{itemize}
			\item {\it Tilahun Abebaw, Department of Mathematics, Addis Ababa University, P. O. BOX 1176, Addis Ababa, Ethiopia, {\color {orange} tilahun.abebaw@aau.edu.et}}
			\item {\it Nega Arega, Department of Mathematics and Statistics, The Namibia University of Science and Technology, Namibia,  {\color{orange}nechere@nust.na}}
			\item {\it Teklemichael Worku Bihonegn, Department of Mathematics, Addis Ababa University, P. O. BOX 1176, Addis Ababa, Ethiopia, {\color {orange}teklemichael.worku@aau.edu.et}}
			\item {\it David Ssevviiri, Department of Mathematics, Makerere University, P. O. BOX 7062, Kampala, Uganda, {\color {orange}david.ssevviiri@mak.ac.ug}}
		\end{itemize}


\addcontentsline{toc}{chapter}{Bibliography}
		\begin{thebibliography}{99}
			\bibitem{agnarsson2020monomial} G. Agnarsson and N. Epstein, On monomial ideals and their socles, Order {\bf 37} (2) (2020), 341--369.
			
			\bibitem{agnarsson2023posets} G. Agnarsson and N. Epstein, On posets, monomial ideals, Gorenstein ideals and their combinatorics,
			arXiv preprint arXiv:2302.10068, (2023).
			
			\bibitem{azizi2007radical} A. Azizi, Radical formula and prime submodules, J. Algebra {\bf 307} (1) (2007), 454--460.
			\bibitem{azizi2009radical} A. Azizi, Radical formula and weakly prime submodules, Glasgow Math. J. {\bf 51} (2) (2009), 405--412.
			\bibitem{bhatt2019prisms} B. Bhatt 
			and P. Scholze, Prisms and prismatic cohomology, Ann. of Math. {\bf 196}(3) (2022), 1135-1275.
			\bibitem{bruns1998cohen} W. Bruns and H. J. Winfried, Cohen-macaulay rings, Cambridge university press, Cambridge (39) (1998).
			\bibitem{eizenbud1995commutative} D. Eisenbud, Commutative algebra with a view towards algebraic geometry, Graduate Texts in Mathematics, Springer-verlag, New York, {\bf 150}, (1995).
			\bibitem{elias2021constructive} J. Elias and M. Rossi, A constructive approach to one-dimensional Gorenstein $k$-algebras, Trans. Amer. Math. Soc. {\bf 374} (7) (2021), 4953--4971.
			\bibitem{Elias2012isomorphism} J. Elias and M. E. Rossi, Isomorphism classes of short Gorenstein local rings via Macaulay’s inverse system, Trans. Amer. Math. Soc. {\bf 364} (9) (2012), 4589–4604.
			\bibitem{emsalem1995inverse} J. Emsalem and A. Iarrobino, Inverse system of a symbolic power, I, J. Algebra {\bf 174} (3) (1995), 1080--1090.
			\bibitem{fulton1997young} W. Fulton, Young tableaux: with applications to representation theory and geometry, Mathematical Society Student Texts, Cambridge University Press, Cambridge, {\bf 35} (1997).
			
			\bibitem{geramita1996waring}A. V. Geramita, Inverse systems of fat points: Waring’s problem, secant varieties of Veronese varieties and parameter spaces for Gorenstein ideals., The Curves Seminar at Queen’s {\bf 40} (13) (1996).
			\bibitem{harbourne2011inverse}B. Harbourne, H. Schenck, and A. Seceleanu, Inverse systems, Gelfand-Tsetlin patterns and the weak
			Lefschetz property, J. Lond. Math. Soc. {\bf 84} (3) (2011), 712–730.
			
			\bibitem{Jenkins1992} J. Jenkins and P. F. Smith, On the prime radical of a module over a commutative ring, Comm. Algebra {\bf 20} (12) (1992), 3593--3602.
			
			\bibitem{cuong2018commutative}E. Juan, Inverse systems of local rings, 
			Comm. Algebra and its Interactions to Algebraic Geometry: VIASM 2013--2014, {\bf 2210} (2018), 119--163.
			\bibitem{kimuli2022characterizations} P. I. Kimuli and D. Ssevviiri, Characterizations of regular modules, Int. Electron. J. Algebra {\bf 33} (2023), 54--76.
			
			\bibitem{kyomuhangi2021generalised}A. Kyomuhangi and D. Ssevviiri, Generalized reduced modules, Rend. Circ. Mat. Palermo(2) {\bf 72} (1) (2023), 421--431.	
			\bibitem{kyomuhangi2020locally}A. Kyomuhangi and D. Ssevviiri, The locally nilradical for modules over commutative rings, Beitr.  Algebra Geom. {\bf 61}(4) (2020), 759--769.
			
			\bibitem{lee2019semiprime} S. C. Lee and R. Varmazyar, Semiprime submodules of a module and related concepts, J. Algebra Appl. {\bf 18} (08) (2019).
			\bibitem{lee2004reduced}T. K. Lee and Y. Zhou, Reduced modules, 
			rings, modules, algebras and abelian groups, Lecture Notes in Pure and Appl. Math. {\bf 236} (2004), 365--377.
			\bibitem{leung1997}K. H. Leung and S. H. Man, On commutative Noetherian rings which satisfy the radical formula, Glasgow Math. J. {\bf 39} (1997), 285--293.
			
			\bibitem{macaulay1994algebraic} F. S. Macaulay, The algebraic theory of modular systems, Cambridge Math. Lib. {\bf 19} (1994).
			\bibitem{mccasland1991radicals} R. L. McCasland and M. E. Moore, On radicals of submodules, Comm. Algebra {\bf 19} (5) (1991), 1327--1341.
			
			\bibitem{pusat1986} D. Pusat-Yilmaz and P. F. Smith, Modules which satisfy the radical formula, Acta. Math. Hungar. {\bf 95} (2002), 155-167.
			\bibitem{rege2008reduced} M. B. Rege and A. M. Buhphang, On reduced modules and rings, Int. Electron. J. Algebra {\bf 3} (2008), 58--74.
			\bibitem{rohrer2019torsion} F. Rohrer Torsion functors, small or large,
			Beitr. Algebra Geom. {\bf 60}(2) (2019), 233--256.
			
			\bibitem{schenzel2018completion} P. Schenzel and A. M. Simon, Completion, {\v{C}}ech and local homology and cohomology, interactions between them, Springer Monographs in Mathematics, Springer, Cham (2018).
			
			\bibitem{sharif2002} H. Sharif, Y. Sharifi and S. Namazi, Rings satisfying the radical formula, Acta. Math. Hungar. {\bf  71} (1996), 103--108.
			
			\bibitem{ssevviiri2022applications} D. Ssevviiri, Applications of reduced and coreduced modules I, Int. Electron. J. Algebra (2023), DOI: 10.24330/ieja.1299587.
			\bibitem{ssevviiri2023applications} D. Ssevviiri,
			Applications of reduced and coreduced modules II, 
			arXiv preprint arXiv:2205.13241, (2023).
			
			\bibitem{david2015} D. Ssevviiri, A relashionship between 2-primal modules and modules that satisfy the radical formula, Int. Electron. J. Algebra {\bf 18} (2015), 34--45.
			
			\bibitem{van2010simplicial}A. Van Tuyl and F. Zanello, Simplicial complexes and Macaulay’s inverse systems, Math. Z. {\bf 265} (2010), 151--160.
			\bibitem{villarreal7300monomial} R. H. Villareal
			Monomial Algebras, Monographs and Textbooks in Pure and Applied Mathematics, Marcel Dekker, New York, {\bf 238} (2001).
			\bibitem{wolff2016survival} A. R. G. Wolff, The survival complex, arXiv preprint arXiv:1602.08998, (2016).
			
			
			\bibitem{yekutieli2021weak} A. Yekutieli, Weak proregularity, derived completion, adic flatness, and prisms, J. Algebra {\bf 583} (2021), 126--152.
		\end{thebibliography}
	\end{document}